%Paper: 9205001
%From: piccoli@tsmi19.sissa.it
%Date: Wed, 06 May 1992 17:03:24 +0100

\font\bigbf=cmbx10 at 16pt
\hfuzz=20pt
\font\medbf=cmbx10 at 13pt
\magnification=\magstep1
\def \n{\noindent}
\def \no{{\bf n}}
\def \v{\vskip 1em}
\def \m{{m^*}}
\def \vs{\vskip 2em}
\def \vsk{\vskip 4em}
\def \R {I\!\!R}

\def \ve{\varepsilon}
\def \Sab{{S^{b,q}_{a,p}}}
\def \S {{\cal S}}

\def \K {{\cal K}}

\baselineskip = 18 truept
\null
\vsk
\vsk
\centerline {\bigbf A Baire Category Approach}
\v
\centerline {\bigbf to the Bang-Bang Property}
\vsk
\centerline {\it Alberto Bressan and Benedetto Piccoli}
\v
\centerline {S.I.S.S.A., Via Beirut 4,}

\centerline{34014 Trieste, Italy}

\vsk
\item {~} {\it Abstract} -  Aim of this paper is to develop a new technique,
based on the Baire category theorem, in order to establish the closure
of reachable sets and the existence of optimal trajectories for control
systems, without the usual convexity assumptions. The bang-bang property
is proved for a new class of ``concave" multifunctions, characterized
by the existence of suitable linear selections. The proofs rely on
Lyapunov's theorem in connection with a Baire category argument.
\vsk
\centerline{Ref. S.I.S.S.A. ~ 70/92/M}
\vfill\eject
\n
{\medbf 1 - Introduction.}
\v
Aim of this paper is to develop a new technique, based on the Baire category
theorem, in order to establish the closure of reachable sets and the
existence of optimal trajectories for control systems, without the usual
convexity assumptions.

Most of our results will be formulated within the framework of differential
inclusions.  Let $F:\R\times\R^n\mapsto 2^{\R^n}$ be a continuous
multifunction with compact convex values and denote by $ext F(t,x)$ the set
of extreme points of $F(t,x)$.  We say that $F$ has the {\sl bang-bang
property} if, for every interval $[a,b]$ and every Caratheodory solution
$x(\cdot)$ of
$$\dot x(t)\in F(t,x(t))\qquad\qquad t\in [a,b],\eqno(1.1)$$
there exists also a solution of
$$\dot y(t)\in extF(t,y(t))\qquad\qquad t\in[a,b]\eqno(1.2)$$
such that
$$y(a)=x(a),\qquad\qquad y(b)=x(b).\eqno(1.3)$$
If $A,B$ are respectively $n\times n$ and $n\times m$ matrices, and $U
\subset\R^m$ is compact convex, the well known Bang-Bang Theorem [8, 10, 15]
implies that the above property holds for the ``linear" multifunction
$$F(t,x)=\big\{A(t)x+B(t)u;~~~~u\in U\big\}\subset\R^n.\eqno(1.4)$$
In the present paper, the bang-bang property is proved for a new
class of ``concave"
multifunctions, characterized by the
existence of suitable linear selections.  The proofs rely on
Lyapunov's theorem in connection with a Baire category argument.
As applications, we obtain some closure theorems for the reachable set of a
differential inclusion with non-convex right hand side, and new existence
results for optimal control problems in Mayer as well as in Bolza form.

Roughly speaking, the Baire category method consists in showing that the
set $\S^{ext F}$ of solutions of (1.2) is the intersection of countably
many relatively open and dense subsets of the
family $\S^F$ of all solutions of (1.1).  Since $\S^F$ is closed, Baire's
theorem thus implies $\S^{extF}\not=\emptyset$.  The effectiveness of such an
argument, in connection with the Cauchy problem for a differential inclusion,
was suggested by Cellina [5] and demonstrated in [4, 9, 19] and in other
papers.
Here, this basic technique will be combined with
Lyapunov's theorem and applied to the two-point boundary
value problem (1.2), (1.3).

The use of a Lyapunov-type theorem, in order to prove existence of optimal
solutions for non-convex control problems, was introduced by
Neustadt [12] and later applied in [1, 13, 16] to a variety of
optimization problems,
always in connection with evolution equations and cost functionals which
are linear w.r.t. the state variable.  In [6], Cellina and Colombo showed
that the linear cost functional can be replaced by one
which is concave w.r.t. the state variable.  Extensions
and applications to partial differential equations have recently appeared
in [7, 14].  We remark that, if a variational problem of the type
considered in [6] is reformulated as a
Mayer problem of optimal control, then the corresponding multifunction
satisfies our concavity assumptions.  The present results can thus be
regarded as a natural extension of the theorem in [6], for optimization
problems which are ``fully concave": in their dynamics
as well as in the cost functional.

\vsk
\n {\medbf 2 - Preliminaries.}
\v
In this paper, $|\cdot|$ is the euclidean norm in $\R^n$,
$B(x,r)$ denotes the open ball centered at $x$ with radius $r$,
while $B(A,\ve)$ denotes the open $\ve$-neighborhood around the set $A$.
We write $\overline A$ and $\overline{co} A$ respectively for the closure
and the closed convex hull of $A$, while $A
\setminus B$ indicates a set-theoretic difference.
The Lebesgue measure of a set $J\subset\R$ is $meas(J)$.
We recall that a subset $A\subseteq S$ is a $G_\delta$ if $A$ is
the intersection of countably many relatively open subsets of $S$.
\v
In the following, $\K_n$ denotes the family of all nonempty compact convex
subsets of $\R^n$, endowed with the Hausdorff metric.  A key technical tool
used in our proofs will be the function
$h:\R^n\times\K_n\mapsto \R\cup\{-\infty\}$, defined by
$$h(y,K)=sup\left\{\left(\int_0^1| f(x)-y|^2\ dx\right)^{1\over 2}
;~~\ f:[0,1]\rightarrow K,\ \int_0^1 f(x)\ dx=y\right\},\eqno(2.1)$$
with the understanding that $h(y,K)=-\infty$ if $y\notin K$.
Observe that $h^2(y,K)$ can be interpreted as the maximum variance among all
random variables supported inside $K$, whose mean value is $y$.
{}From the above definition, it is clear that
$$h(\xi+y, ~\xi+K)=h(y,K),\qquad ~~~h(\lambda y,~\lambda K)=\lambda h(y,K),
\qquad \forall \xi\in\R^n,~\lambda>0.\eqno(2.2)$$
\v
For the basic theory of multifunctions and differential inclusions we refer
to [1].  Given a solution $x(\cdot)$ of (1.1), following [3] we define its
{\it likelihood} as
$$ L(x)=\left(\int_a^b h^2(\dot x(t),F(t,x(t))dt\right)^{1\over 2}=
\|h(\dot x,~F(\cdot,x)\|_{{\cal L}^2}.\eqno(2.3)$$
\v
The following results were proved in [3]:
\v
\n{\bf Lemma 1.} {\it For every $y,K$, one has $h(y,K)\leq r(K)$,
where $r(K)$ is the radius of the smallest ball containing $K$ (i.e.,
the Cebyshev radius).  Moreover, $h(y,K)=0$ iff $y\in ext K$. Therefore,
a solution $x(\cdot)$ of (1.1) satisfies also (1.2) iff $L(x)=0$.}
\v
\n{\bf Lemma 2.} {\it The map $(y,K)\mapsto h(y,K)$ is upper semicontinuous
in both variables and concave w.r.t. $y$.  The map $x(\cdot)\mapsto
L(x)$ is upper semicontinuous on the set of solutions of (1.1), endowed
with the ${\cal C}^0$ norm.}
\vsk
\n{\medbf 3 - The main results.}
\v
In the following, we denote by $\Sab$ the set of
all Caratheodory solutions of the two-point boundary value problem
$$\dot x(t)\in F(t,x(t)),\qquad\qquad x(a)=p,~\qquad x(b)=q.\eqno(3.1)$$
\v
\noindent{\bf Theorem 1.}
{\it Let $F:\R\times\R^n\mapsto 2^{\R^n}$ be a continuous multifunction
with compact, convex values.
The following conditions are equivalent:
\item {\rm (1)} For every interval $[a,b]$ and every $p,q \in \R^n$,
if $\Sab\not=\emptyset$, then
the set of solutions of
$$\dot y(t)\in  extF(t,y(t)),\qquad\qquad x(a)=p,~~~x(b)=q\eqno(3.2)$$
is a dense $G_\delta$ in $\Sab$.
\item {\rm (2)} $F$ has the bang-bang property
\item {\rm (3)} For every interval $[a,b]$ and every $p,q \in \R^n$, if
$\Sab \neq \emptyset$, then for every $\ve >0$ there exists a solution
$x(\cdot)$ of (3.1) such that
$$L^2(x)\doteq\int_a^{b}\ h^2\big(\dot x (t),F(t,x(t))\big)~dt~ <~\ve.
\eqno(3.3)$$

\n Proof.}
(1) $\Rightarrow$ (2)~~ If $\Sab\neq\emptyset$, then by (1) the
set of solutions of (3.2), being dense, is nonempty. Hence (2) holds.
\v
\n (2) $\Rightarrow$ (3)~~ If $\Sab \neq \emptyset$ then by (2) there
exists a solution $y(\cdot)$ of (3.2). This implies (3), because by Lemma
1
$$\int_a^b\ h^2(\dot y (t),F(t,y(t)))\ dt = 0<\ve.$$
\v
\n(3) $\Rightarrow$ (1)~~Consider the sets $A_m\doteq\big\{x\in\Sab~;~~ L(x)<{
1 \over m}\big\}$.  By Lemma 2, $L$ is upper semicontinuous, hence each
$A_m$ is open.
Now fix any $x(\cdot)\in\Sab$,
$\ve >0$.  Define
$$\Omega\doteq\big\{(t,z);~~~t\in [a,b],~~|z-x(t)|\leq\ve\big\}$$ and choose a
constant $M$ so large that
$$F(t,x)\subseteq \overline B(0,M)\qquad\qquad\forall (t,x)
\in\Omega.\eqno(3.4)$$
Split the interval $[a,b]$ into $k$ equal subintervals $J_i=[t_{i-1},~t_i]
$, inserting the points $t_i\doteq a+(i/k)(b-a)$, choosing $k$
so large that $ 2M(b-a)/k\leq\ve$.

By the assumption (3), for each $i$ there exists a solution $y_i:[t_{i-1},
{}~t_i]\mapsto\R^n$ of the two-point boundary value problem
$$\dot y(t)\in F(t,y(t)),\qquad~~~~y(t_{i-1})=x(t_{i-1}),~~~y(t_i)=
x(t_i),\eqno(3.5)$$
with
$$L^2(y_i)=\int_{t_{i-1}}^{t_i} h^2\big(\dot y_i(t),~F(t,y_i(t))\big)~dt
<{1\over m^2k}.\eqno(3.6)$$
Define $y(\cdot)$ as the solution of (3.1) whose restriction to
each $J_i$ coincides with $y_i$.
Given any $t\in [a,b]$, if, say, t$\in J_i$, then (3.4), (3.5) imply
$$|y(t)-x (t)|\leq\int_{t_i}^t
|\dot y_i(s)-\dot x(s)|~ ds\leq {2M(b-a)\over k}
\leq\ve.$$ Hence $\|y-x\|_{{\cal C}^0}\leq\ve$.  Moreover, $y\in A_m$
because
$$L^2(y)=\sum_{i=1}^k \int_{t_{i-1}}^{t_i} h^2\big(\dot y_i(t),~F(t,y_i(t))
\big)~dt<{k\over m^2k}={1\over m^2}.$$
Since $x(\cdot)$ and $\ve>0$ were arbitrary, this proves that each
$A_m$ is dense in $\Sab$.  By Baire's theorem,
it follows that
$A =\bigcap_m A_m$ is a $G_\delta$ dense subset of $\Sab$.
If $y\in A$, then $L(y)=0$ and hence $\dot y(t)\in extF(t,y(t))$ almost
everywhere.
\vs
In the previous theorem, the implication (3) $\Rightarrow$ (1) determines
the strength of the category method.  In order to prove that ``most"
solutions of (3.1) actually solve (3.2) as well, it suffices to show
(for every $a,b,p,q$) the
existence of some solution of (3.1) with arbitrarily small likelihood.
Roughly speaking, this requires the construction of some solution $y$
of (3.1) whose derivative remains close to the extreme points of $F(t,y)$
during most of the time.

In practice, the condition (3) may often
be easier to verify.  We now show that this is indeed the case, if the
multifunction $F$ satisfies suitable concavity conditions.
\vs
\n{\bf Theorem 2.} {\it Let $F:\R\times\R^n\mapsto 2^{\R^n}$ be a Hausdorff
continuous multifunction with compact, convex values.
Assume that:
\v
\item {\rm (C1)}  For each $(t,x)$ and every $y\in F(t,x)$,
there exists
a linear function $z\mapsto Az+c$ satisfying
$$y =Ax+c,\qquad\qquad
Az+c\in F(t,z)\qquad\forall z\in \overline B(x,\rho(t,x)),\eqno(3.7)$$
where the radius $\rho=\rho(t,x)$ remains uniformly positive on compact sets.
\v
\item {\rm (C2)} For each $(t,x)$, every $y\in F(t,x)$ and $\ve>0$,
there exist $\delta>0$ and $n+1$ linear functions $z\mapsto A'z+c_i$,
$~~i=0,\ldots,n$, such that
$$y\in\overline{co}\big\{A'x+c_0~,~\ldots~,~A'x+c_n\big\},\eqno(3.8)$$
$$h(A'x+c_i,~F(t,x))\leq\ve\quad\qquad\forall i,\eqno(3.9)$$
$$ A' z+c_i\in F(t,z)\qquad\quad\forall z\in\overline B(x,\delta),~\forall
i.\eqno(3.10)$$

\n Then F has the bang-bang property.}
\vs
We refer to (C1), (C2) as {\it concavity conditions} because they require, for
each point $(t,x,y)$ of the graph of $F$, the existence of suitable linear
(non-homogeneous) selections.  A similar property is shared by the
epigraph of a concave scalar function, which admits global
linear selections through each of its points.
\v
\n {\it Proof of Theorem 2.}
{}~ We will prove that $F$ has property (3) stated in Theorem 1.
Let $x^*(\cdot)$ be a solution of (3.1), for some interval $[a,b]$ and some
points $p,q\in\R^n$.  Let any $\ve >0$ be given, and define
$$\eta=\inf_{t\in [a,b]}~ \rho(t,x^*(t)),\qquad\qquad V=\big\{(t,z);~~
t\in [a,b],~~|z-x^*(t)|\leq\eta\big\}.\eqno(3.11)$$
By assumption, $\eta>0$.  Choose $M$ so large that
$$F(t,z)\subseteq\overline B(0,M)\qquad\qquad\forall (t,z)\in V.\eqno(3.12)$$
By Lemma 1, this implies
$$h(y,~F(t,z))\leq M\qquad\qquad\forall y,~~\forall (t,z)\in V.\eqno(3.13)
$$
\v
\n{\bf 1.}  As a first step, we construct measurable, bounded functions
$A,~c$, such that
$$\dot x^*(t)=A(t)x(t)+c(t)\qquad\quad\hbox{for a.e.~}t\in[a,b],
\eqno(3.14)$$
$$A(t)z+c(t)\in F(t,z)\qquad\qquad\forall z\in\overline B(x^*(t),\eta).
\eqno(3.15)$$
Since $\dot x^*(\cdot)$ is measurable, by Lusin's theorem
there exists a sequence of disjoint compact sets $(J_\nu)_{\nu\geq 1}$
such that
$$meas\Big([a,b]\setminus\bigcup_{\nu\geq 1}J_\nu\Big)=0\eqno(3.16)$$
and such that the restriction of $\dot x^*$ to each $J_\nu$ is continuous.
Define the multifunction ~$G:[a,b]\mapsto 2^{\R^{n\times n}
\times\R^n}$ by setting
$$G(t)\doteq\big\{ (A,c):\quad \dot x^*(t)=Ax+c,~~~Az+c\in F(t,z)~~~
\forall z\in\overline B(x^*(t),\eta)\big\}.$$
Because of (C1) and of the choice of $\eta$,~ $G(t)\not=\emptyset$ for
a.e. $t$.
One easily checks that the restriction of $G$ to each $J_\nu$ has closed
graph, because of the continuity of $\dot x^*,~x^*$ and $F$.
Hence, $G$ is a measurable multifunction on $[a,b]$ with closed, nonempty
values.  By [11], it admits a measurable selection $t\mapsto \big(A(t),~
c(t)\big)$, which clearly satisfies (3.14), (3.15).
Observe that the matrices $A(t)$ and the vectors $c(t)$ must be uniformly
bounded, because of (3.15), (3.12).
\v
\n{\bf 2.} As a second step, we construct measurable functions
$A',c_0,\ldots c_n,\theta_0,\ldots,\theta_n,\delta$, such that, for
almost every $t\in [a,b]$, the following holds:
$$\delta(t)>0,\qquad\theta_i(t)\in [0,1],\qquad\sum_{i=0}^n\theta_i(t)=1,
\eqno(3.17)$$
$$\dot x^*(t)=A'(t)x^*(t)+\sum_{i=0}^n \theta_i(t)c_i(t),\eqno(3.18)$$
$$A'(t)z+c_i(t)\in F(t,z),\qquad h^2\big(A'(t)z+c_i(t),~F(t,z)\big)
\leq\ve\qquad\forall i,~~\forall z\in\overline B\big(x^*(t),\delta(t)
\big).\eqno(3.19)$$

By Lemma 2, $h$ is upper semicontinuous, hence there exists a
nonincreasing sequence $(h_m)_{m\geq 1}$ of continuous
functions such that
$$h(y,K)=\inf_{m\geq 1}~ h_m(y,K)\qquad\qquad \forall (y,K)\in \R^n
\times\K_n.\eqno(3.20)$$
For each $m\geq 1$, define the multifunction
$$\eqalign{H_m(t)\doteq \bigg\{~&\big( A',c_0,\ldots,c_n,\theta_0,\ldots,
\theta_n,\delta\big);\qquad\quad \delta={1\over m},\qquad \theta_i\in [0,1]
,\cr &\qquad\dot x^*(t)=\sum_{i=0}^n\theta_i(A'x(t)+c_i),\qquad
\sum_{i=0}^n\theta_i=1,\cr
&A'z+c_i\in F(t,z),\qquad h^2_m\big(A'z+c_i,~F(t,z)\big)\leq\ve
\qquad\forall z\in\overline B\Big(x(t),{1\over m}\Big)\bigg\}.\cr}$$
If $(J_\nu)_{\nu\geq 1}$ is the same sequence of compact sets considered
at (3.16), the continuity of $\dot x^*$, $x^*$, $h_m$ and $F$
implies that the restriction of $H_m$ to each $J_\nu$ has closed graph,
with uniformly bounded, possibly empty values.

Defining $I_m=\{t;~H_m(t)\not=\emptyset\}$, it is clear that on each
$I_m$ the multifunction $H_m$ is measurable with closed, nonempty values.
By [11], it admits a measurable selection, say $t\mapsto \Phi_m(t)$.
By (C2), (3.20) and the continuity of each $h_m$, for every $\nu$ we have
$\bigcup_{m\geq 1} I_m\supseteq J_\nu$.
Therefore, the selection
$$\big(A'(t),c_0(t),\ldots,c_n(t),\theta_0(t),\ldots,\theta_n(t),
\delta(t)\big)\doteq\Phi_m(t)\qquad\hbox{iff}\quad t\in I_m\setminus
\bigcup_{\ell<m}I_\ell$$
is measurable and defined for a.e. $t\in [a,b]$. By construction, the
conditions (3.17)-(3.19) hold.
\v
\n{\bf 3.} We can now complete the proof of the theorem.
Since $\delta(\cdot)$ is measurable and positive, there exists an integer
$\m$ such that
$${1\over {\m}}< \eta,\qquad
meas\big(J'_{\m}\big)\geq (b-a)-\varepsilon,\eqno(3.21)$$
where
$$J'_m\doteq\Big\{ t\in[a,b];~~\delta(t)\geq{1\over m},~~~|A'(t)|,
|c_i(t)|\leq m\Big\}.\eqno(3.22)$$
Split $[a,b]$ into $k$ equal subintervals $I_j=[t_{j-1},~t_j]$, inserting
the points $t_j\doteq a+(j/k)(b-a)$ and choosing $k$ so large that:
$$2M{b-a\over k}< {1\over m^*}.\eqno(3.23)$$
Using the selections $A,~c$~ and~ $A',~c_i,~\theta_i~$
constructed in the previous steps, define:
$$A^*(t) =\cases{A(t)\quad&if\qquad $t\notin J_{\m}',$\cr &\cr
                A'(t)\quad&if\qquad $t\in J_{\m}',$\cr}$$
$$f(t) =\cases{c(t)\qquad &if\qquad $t\notin J_{\m}'$,\cr &\cr
\sum_{i=0}^n\theta_i(t)c_i(t)\quad&if\qquad $t\in J_{\m}'$.\cr}$$
By (3.14), (3.18),~ $\dot x^*(t)=A^*(t)x(t)+f(t)$.
Calling $W(\cdot,\cdot)$
the matrix fundamental solution of the bounded linear system
$\dot v=A^*(t)v$, we thus have the representation
$$x^*(t)=W(t,t_{j-1})x^*(t_{j-1})+\int_{t_{j-1}}^tW(t,s)f(s)~ds,
\qquad ~~t\in [t_{j-1},~t_j].$$
Applying Lyapunov's theorem on each interval $I_j$,
for every $j$ we obtain a measurable partition $\big\{ I_{j,0},
\ldots,I_{j,n}\big\}$ of $I_j$ and an absolutely continuous function
$w_j$ satisfying the two-point boundary value problem
$$\dot w_j(t)=\cases{ A^*(t)w_j(t)+c(t)\qquad &if\qquad $t\in I_j\setminus
J_{m^*}'$,\cr &\cr A^*(t)w_j(t)+c_\ell(t)\qquad &if\qquad $t\in I_{j,\ell}
\cap J'_{m^*},~~\ell=0,\ldots,n$,\cr}$$
$$w_j(t_j)=x^*(t_j),\qquad w_j(t_{j+1})=x^*(t_{j+1}).$$

We claim that
$$|w_j(t)-x^*(t)|\leq{1\over {\m}}\qquad\qquad\forall t\in I_j.\eqno(3.24)$$
If not, there would exist a first time $\tau\in I_j$ such that
$$|w_j(\tau)-x^*(\tau)|={1\over {\m}}.\eqno(3.25)$$
Recalling (3.15), (3.19) and using (3.12) and (3.23), we then have:
$$|w_j(\tau)-x^*(\tau)|\leq\int_{t_{j-1}}^\tau|\dot w_j(t)
-\dot x^*(t)|\ dt\leq 2M{b-a\over k}<{1\over \m},$$
a contradiction with (3.25). This proves (3.24). In particular,
by (3.15), (3.19) we conclude that $\dot w_j(t)\in F(t, w_j(t))$
for a.e. $t\in I_j$.

Now consider the solution $w(\cdot)$ of (3.1) whose restriction to each
$I_j$ coincides with $w_j$.
Recalling (3.13), (3.19), (3.21), its likelihood is computed by
$$\eqalign{L^2(w) &= \int_{J'_{m^*}} h^2(\dot w(t),F(t,w(t)))~dt+
\int_{[a,b]\setminus J'_{m^*}} h^2(\dot w(t),F(t,w(t)))~dt\cr
&\leq\ve\cdot meas\big(J_{m^*}'\big)+M^2\cdot meas\big([a,b]\setminus
J_{m^*}'\big)
\leq\ve ((b-a)+M^2).\cr}$$
Since $\ve$ was arbitrary, this establishes the property
(3) in Theorem 1, which is equivalent to the bang-bang property.
\v
\n{\it Remark 1.}  The previous theorems, with the same proofs, remain
valid if $F$ is defined on some open set $\Omega\subset \R\times\R^n$.
\vsk
\n{\medbf 4 - Examples of concave multifunctions.}
\v
Aim of this section is to exhibit some classes of multifunctions
which satisfy the concavity properties (C1), (C2) stated in Theorem 2.
\v
\n{\bf Proposition 1.}  {\it Let $\varphi:\R\times\R^n\mapsto ]0,\infty[$
be a continuous function, with $x\mapsto \varphi(t,x)$ convex for every
$t$.  Let $U\subset\R^n$ be compact, convex, containing the origin.
Then the multifunction
$$F(t,x)=\varphi(t,x)U\eqno(4.1)$$
has the bang-bang property.
\v
\n Proof.}  In order to apply Theorem 2, we first verify the concavity
condition (C1).  Fix any $(t,x)$ and any $y=\varphi (t,x)u\in
F(t,x)$.  Since $\varphi$ is continuous and strictly positive and its
subdifferential $\partial_x \varphi$ ~w.r.t. $x$ is uniformly bounded on
compact sets, we have
$$\varphi(t,x)+\xi\cdot(z-x)\geq 0\qquad~\forall\xi\in\partial_x\varphi
(t,x),\quad\forall z\in\overline B(x,\rho(t,x)),\eqno(4.2)$$
for some function $\rho=\rho(t,x)$ uniformly positive on compact sets.
Choose any vector $\xi\in\partial_x
\varphi(t,x)$ and define the linear map
$$z\mapsto Az+c\doteq (\xi\cdot z)u+\big[ \varphi(t,x)u-(\xi\cdot x)u\big].
$$
If $|z-x|\leq \rho(t,x)$, we need to show the existence of some
$\omega\in U$ such that
$$(\xi\cdot z)u+\big[\varphi(t,x)u-(\xi\cdot x) u\big]=\varphi(t,z)\omega.
\eqno(4.3)$$
{}From (4.2) and the convexity of $\varphi$ it follows
$$\omega=
{\varphi(t,x)+\xi\cdot (z-x)\over \varphi (t,z)} u=\alpha u\eqno(4.4)$$
for some $\alpha\in [0,1]$.  The assumptions on $U$ thus imply $\omega\in
U$.
\v
\n We now turn to the condition (C2).  Let
$(t,x)$,  $y=\varphi(t,x)u\in F(t,x)$ and $\ve>0$ be given.  We can assume
$$u=\sum_{i=1}^\nu \theta_i u_i,\qquad~~\theta_i\in (0,1],\qquad
\sum_{i=1}^\nu\theta_i=1\eqno(4.5)$$
for some $\nu\in\{1,\ldots,n+1\}$,~~$u_i\in ext U$.
Select $\xi\in\partial_x\varphi(t,x)$ as before and define
$$u_i'=u_i+\ve'(u-u_i),\eqno(4.6)$$
choosing $\ve'\in (0,1)$ so small that
$$h(u_i',~U)<{\ve\over \varphi(t,x)}\qquad~~\forall i.\eqno(4.7)$$
This is possible because $h(u_i,U)=0$ and $h$ is upper semicontinuous.
Then define
$$A'z=(\xi\cdot z)u,\qquad c_i=\varphi(t,x)u_i'-(\xi\cdot x)u.$$
Recalling (2.2)$_2$, (4.7) yields
$$h\big(A'x+c_i,~F(t,x)\big)=h\big(\varphi(t,x)u_i',~\varphi(t,x)U\big)
<\ve.$$
Hence, for $z$ in a small neighborhood of $x$, the upper semicontinuity
of $h$ implies
$$h\big(A'z+c_i,~\varphi(t,z)U\big)<\ve\qquad\qquad\forall i.$$
Moreover, by (4.5), (4.6),
$$y=\varphi(t,x)u=\sum_{i=1}^\nu\theta_i\varphi(t,x)u_i'\in
\overline{co}\big\{A'x+c_i,~~~i=1,\ldots,\nu\big\}.$$
It remains to prove that
$A'z+c_i\in \varphi(t,z)U$
for $|z-x|$ small enough.  For each fixed $i$, define
the subspace $E_i\doteq span\{u,~u_i\}$.
\v
If $E_i$ has dimension 2,
consider the triangle $\Delta=\overline{co}\{0,~u,~u_i\}$ and call
$\no_1,~\no_2,~\no_3$ the unit vectors in $E_i$ which are outer normals
to the sides $u_i-u$, $u_i$, $u$, respectively.
Observe that
$$\varphi(t,z)\overline{co}\{0,~u,u_i\}=\big\{y\in E_i;~~\no_1\cdot y\leq
\varphi(t,z)(\no_1\cdot u),~~\no_2\cdot y\leq 0,~~\no_3\cdot y\leq 0
\big\}\subseteq F(t,z).\eqno(4.8)$$
Since $u_i'$ is a strict convex combination of $u$ and $u_i$, one has
$\no_2\cdot u_i'<0$, ~$\no_3\cdot u_i'<0$.  By continuity, for $z$
sufficiently close to $x$ we still have
$$\no_j\cdot \big[ \varphi(t,x)u_i'+\xi\cdot(z-x)u\big]=
\no_j\cdot \big[A'z+c_i\big]<0,\qquad j=2,3.\eqno(4.9)$$
Moreover, since $\no_1\cdot u=\no_1\cdot u_i'>0$, the convexity
of $\varphi$ implies
$$\eqalign{\varphi(t,z)(\no_1\cdot u)&\geq\big[\varphi(t,x)+\xi\cdot(z-x)
\big](\no_1\cdot u)\cr
&=\no_1\cdot\big[\varphi(t,x)u_i'+\big(\xi\cdot (z-x)\big)u\big]
\cr &=\no_1\cdot[A'z+c_i].\cr}\eqno(4.10)$$
By (4.8), the inequalities (4.9), (4.10) together imply $A'z+c_i\in F(t,z)$.
\v
Finally, consider the case where $E_i$ has dimension $\leq 1$.
Then, either $u=u_i$, hence $\nu=1$ and $u\in ext F(t,x)$.  In this case,
the same argument as in (4.3), (4.4) can be used.  Or else $u_i'$ lies
in the relative interior
of the segment $S\doteq\overline{co}\{0,u,u_i\}$. In this case, the map
$z\mapsto A'z+c_i$ takes values inside $E_i$, with
$\varphi^{-1}(t,x)\big[A'x+c_i\big]\in rel~intS$.
By continuity, $\varphi^{-1}(t,z)\big[A'z+c_i\big]\in S~$
for $|z-x|$ small enough.  This proves again that
$A'z+c_i\in F(t,z)$ for every $z$ in a neighborhood of $x$.

An application of Theorem 2 now yields the desired conclusion.
\vs
The next application is concerned with a multifunction $F$ whose values are
polytopes, with variable shape but constant number of vertices.
More precisely, we assume that $F(t,x)$ admits the representations
$$F(t,x)=\overline{co}\big\{y_1(t,x),\ldots,y_N(t,x)\big\},\eqno(4.11)$$
where $y_1,\ldots,y_N$ are the (distinct) vertices of $F(t,x)$, as well
as
$$F(t,x)=\Big\{y\in \R^n~;~~~w_j(t)\cdot y\leq\psi_j(t,x)=\max_{\omega\in
F(t,x)}w_j(t)\cdot\omega,~~~~j=1,..,k\big\}.\eqno(4.12).$$
On the product set of indices $\{1,\ldots,N\}\times\{1,\ldots,k\}$,
we consider the ``incidence" relation
$$i\sim j\qquad\hbox{iff}\qquad w_j(t)\cdot y_i(t,x)=\psi_j(t,x).\eqno(4.13)$$
\vs
\n{\bf Proposition 2.}  {\it Let $F$ be a multifunction admitting
the representations (4.11), (4.12).  Assume that
\v
\item{(i)} $w_j:\R\mapsto \R^n$,~ $\psi_j:\R\times\R^n\mapsto \R$
are continuous functions, $|w_j(t)|\equiv 1$,
each map $x\mapsto \psi_j(t,x)$ is convex.
\v
\item{(ii)}  The relation $\sim$ defined at (4.13) is independent of
$(t,x)$.
\v
Then $F$ has the bang-bang property.
\v
Proof.}  For each $i\in\{1,\ldots,N\}$,
consider the set of indices
$$J_i=\big\{j;~~~w_j(t)\cdot y_i(t,x)=\psi_j(t,x)\big\}.\eqno(4.14)$$
By (ii), this set does not depend on $t,x$.

We begin by checking the condition (C1) in Theorem 2.  First, assume
$y\in extF(t,x)$, say $y=y_i(t,x)$.
In this case we can choose $n$ independent vectors, say
$w_{j_1}(t),\ldots,w_{j_n}(t)$, with $j_1,\ldots,j_n\in J_i$.
Define the dual vectors $w^*_{j_\ell}$, requiring that
$$w^*_{j_\ell}\cdot w_{j_m}=\cases{1\qquad &if\qquad $\ell=m,$\cr
0\qquad &if\qquad $\ell\not= m$.\cr}\eqno(4.15)$$
By convexity, each function $z\mapsto\psi_j(t,z)$ is differentiable almost
everywhere.  Therefore, there exists a sequence of points $x_\nu
\to x$ such that the gradients $\nabla_x\psi(t,x_\nu)$ exist for each
$\nu$, together with the limits
$$\lim_{\nu\to\infty}\nabla_x\psi_j(t,x_\nu)=\xi_j\in
\partial_x\psi_j(t,x)\qquad~~~\forall j.\eqno(4.16)$$
Now define
$$A_i z=\sum_{\ell=1}^n(\xi_{j_\ell}\cdot z)w^*_{j_\ell},\qquad
c_i=y_i(t,x)-\sum_{\ell=1}^n(\xi_{j_\ell}\cdot x)w_{j_\ell}^*.\eqno(4.17)$$
Clearly, $A_ix+c_i=y_i$.  Using the representation (4.12) we now check that
$$A_iz+c_i\in F(t,z)\qquad\quad\forall z\in\overline B(x,\rho(t,x))\eqno
(4.18)$$ for some $\rho=\rho(t,x)$ uniformly positive on compact sets.

If $j\in J_i$, then there exist unique coefficients
$\alpha_\ell$ such that $w_j(t)=\sum \alpha_\ell w_{j_\ell}(t)$.  The
assumption (ii) together with (4.16) now implies
$$\psi_j(t,z)=w_j(t)\cdot y_i(t,z)=\sum_{\ell=1}^n\alpha_\ell
\psi_{j_\ell}(t,z)\qquad~~~\forall z,$$
$$\sum_{\ell=1}^n\alpha_\ell\xi_{j_\ell}\in\partial_x\psi_j(t,x).
\eqno(4.19)$$
{}From (4.19) and the convexity of $\psi_j$ it follows
$$\eqalign{ w_j(t)\cdot \big[A_iz+c_i\big] &=
\sum_{\ell=1}^n\alpha_\ell w_{j_\ell}(t)\cdot\left[ y_i(t,x)+
\sum_{h=1}^n\big(\xi_{j_h}\cdot(z-x)\big)w^*_{j_h}\right]\cr
&=\psi_j(t,x)+\Big(\sum_{\ell=1}^n\alpha_\ell\xi_{j_\ell}
\Big)\cdot(z-x)\leq\psi_j(t,z).\cr}$$
On the other hand, if $j\notin J_i$ then
$w_j(t)\cdot\big[A_ix+c_i\big]<\psi_j(t,x).$
Hence, by continuity
we still have $$w_j(t)\cdot\big[ A_i z+c_i\big] <\psi_j(t,z)
\qquad~~\forall z\in\overline B\big(x,\rho(t,x)\big).$$
By the assumption (ii), the continuity of the functions $w_j,\psi_j$
and the local boundedness of the subgradients $\partial _x\psi_j$,
it follows that $\rho$ can be taken uniformly positive on bounded sets.
This proves (C1) in the case $y\in extF(t,x)$.

When $y$ is an arbitrary element in $F(t,x)$, there exist extreme points
$y_i$ and coefficients $\theta_i\in [0,1]$ such that
$$y=\sum_{i=1}^N\theta_i y_i(t,x),\qquad\quad\sum_{i=1}^N
\theta_i=1.$$
If $A_i, c_i$  are the matrices and vectors defined at (4.17), then
the convex combinations $A=\sum \theta_iA_i,$ ~$c=\sum\theta_ic_i$ satisfy
$$Ax+c=y,\qquad\quad Az+c\in F(t,z)\qquad\forall z\in\overline B(x,\rho(t,x
)).$$
\v
Next, consider the condition (C2).
Let $(t,x)$,  $~y\in F(t,x)$,~ $\ve>0$ be given.  Write $y$ as a convex
combination of points $y_1,\ldots,y_\nu\in extF(t,x)$, say
$$y=\sum_{i=1}^\nu \theta_i y_i\qquad \theta_i \in (0,1]\qquad
\sum_{i=1}^\nu\theta_i =1,$$
and define
$$y_i'=y_i+\ve'(y-y_i),$$
choosing $\ve'\in (0,1]$ so small that
$$h\big(y_i',~F(t,x)\big)<\ve\qquad\quad\forall i.\eqno(4.20)$$
Consider the vector space
$$E\doteq \hbox{span}\big\{ w_j(t);~~~~w_j(t)\cdot y=\psi_j(t,x)\big\}.$$
Choose a basis $\{ w_{j_1},\ldots, w_{j_\mu}\}$ of $E$ and define the dual
basis $\{w_{j_1}^*,\ldots,w_{j_\mu}^*\}$ as in (4.15). Select
vectors $\xi_j\in\partial_x\psi_j(t,x)$ as in (4.16) and define
$$A'z=\sum_{\ell=1}^\mu (\xi_{j_\ell}\cdot z)w^*_{j_\ell},
\qquad~~~c_i=y_i'-\sum_{\ell=1}^\mu(\xi_{j_\ell}\cdot x)w^*_{j_\ell}.
\eqno(4.21)$$
\v
The above definitions imply
$$y=\sum_{i=1}^\nu y_i'\in\overline{co}\big\{A'x+c_i;~~~i=1,\ldots,\mu\big\}.$$
Moreover, by (4.20) and the upper semicontinuity of $h$, for $|z-x|$ small
enough we have
$$h\big( A'z+c_i,~F(t,z)\big)<\ve$$
Using the representation (4.12), we
now prove that $A'z+c_i\in F(t,z)$.
If $w_j(t)\in E$, then
$$w_j(t)\cdot y=w_j(t)\cdot y_i'=\psi_j(t,x)\qquad\qquad\forall i=1,
\ldots,\nu.$$
Moreover, there exist coefficients $\alpha_\ell$ such that
$w_j(t)=\sum \alpha_\ell w_{j_\ell}(t)$.  The assumption (ii) together with
(4.16) now implies
$$\psi_j(t,z)=w_j(t)\cdot y_i(t,z)=\sum_{\ell=1}^\mu\alpha_\ell
\psi_{j_\ell}(t,z)\qquad\quad\forall z,\eqno(4.22)$$
$$\sum _{\ell=1}^\mu \alpha_\ell \xi_{j_\ell}\in\partial_x\psi_j(t,x).\eqno
(4.23)$$
{}From (4.22), (4.23) and the convexity of $\psi_j$ it follows
$$\eqalign{w_j(t)\cdot\big[A'z+c_i\big]&=\sum_{\ell=1}^\mu\alpha_\ell
w_{j_\ell}(t)\cdot\left[ y_i'+\sum_{h=1}^\mu\big(\xi_{j_h}\cdot (z-x)\big)
w_{j_h}^*\right]\cr
&=\psi_j(t,x)+\Big(\sum_{\ell=1}^\mu\alpha_\ell\xi_{j_\ell}\Big)\cdot(z-x)
\cr&\leq\psi_j(t,z).\cr}$$
On the other hand, if $w_j(t)\notin E$, then
$$w_j(t)\cdot [A'x+c_i]=w_j(t)\cdot y_i'<\psi_j(t,x).$$
By continuity, for $|z-x|$ sufficiently small we still have
$$w_j(t)\cdot [A'z+c_i]<\psi_j(t,z).$$
This completes the proof of condition (C2). An application of Theorem 2 now
yields the desired result.
\v
\n{\it Remark 2.}  Assume that $A$, $b$ are a $n\times n$ matrix
and a $n$-vector, depending continuously on $t$, and that $F$ is a continuous,
compact convex valued multifunction satisfying the concavity conditions
(C1), (C2).  Then the multifunction
$$G(t,x)=A(t)x+b(t)+F(t,x)$$
satisfies all assumptions in Theorem 2 as well.  In particular, from
Proposition 1 it follows that the bang-bang property holds for a control
system of the form
$$\dot x=A(t)x+b(t)+\varphi(t,x)u,\qquad~~~u(t)\in U,$$
with $U$ compact, convex, containing the origin and $\varphi>0$ convex w.r.t.
$x$.
\vsk
\n{\medbf 5 - A nonconvex optimal control problem.}
\v
This section is concerned with an application of Theorem 2 to an optimal
control problem of Bolza.  The analysis will clarify the connections
between the concavity conditions (C1), (C2) and the assumptions
made in [6, 12].
Given the linear control system on $\R^n$
$$\dot x(t)=A(t)x+f(t,u(t))\qquad u(t)\in U\eqno(5.1)$$
with initial and terminal constraints
$$x(0)=\bar x\qquad (T,x(T))\in S,\eqno(5.2)$$
consider the minimization problem:
$$\min\ \int_0^T \alpha(t,x(t))+\beta(t,u(t))\ dt.\eqno(5.3)$$
\v
\n{\bf Theorem 3.}  {\it Let the functions $A,f,\alpha,\beta$
be continuous, with $\alpha$ concave w.r.t. $x$.
Assume that the control set $U\subseteq\R^m$ is compact and that the
terminal set $S$ is closed and contained in $[0,T_0]\times\R^n$, for some
$T_0$.  If some solution of (5.1), (5.2) exists, then the minimization
problem (5.3) admits an optimal solution.
\v
\n Proof.}  We begin by adding an extra
variable $x_0$, writing the problem in Mayer form:
$$\min x_0(T)\eqno(5.4)$$
$$\left\{\eqalign{\dot x(t)&=A(t)x(t)+f(t,u(t))\cr
\dot x_0(t)&=\alpha(t,x(t))+\beta(t,u(t))\cr}\right.\qquad
\quad u(t)\in U, \eqno(5.5)$$
$$(x,x_0)(0)=(\bar x,0)\qquad (T,x(T))\in S.\eqno(5.6)$$
The continuity of $A$, $f$ and the compactness of $U$ imply that all
trajectories of the differential inclusion
$$\dot x(t)\in A(t)x(t)+\overline{co}\big\{f(t,u);~~u\in U\big\},
\qquad x(0)=\bar x,\qquad t\in [0,T_0],\eqno(5.7)$$
are contained in some bounded open set $\Omega\subset\R\times\R^n$.
Define the constant
$$M\doteq 1+\sup\big\{\alpha(t,x)+\beta(t,u);\qquad (t,x)\in\Omega,~~u\in U
\big\}\eqno(5.8)$$
and the multifunction (independent of $x_0$)
$$F(t,x,x_0)\doteq \overline{co}\big\{(y,y_0);~~~~y=A(t)x+f(t,u),~~~~
\alpha(t,x)+\beta(t,u)\leq y_0\leq M~~\hbox{for some}~~u\in U\big\}.
\eqno(5.9)$$
Observe that $F$ admits the representation
$$\eqalign{F(t,x)=\bigg\{ (y,y_0);~~~y=A&(t)x+\sum_{i=0}^{n+1}\theta_if(t,
u_i),~~y_0=\Big[\alpha(t,x)+\sum_{i=0}^{n+1}\theta_i\beta(t,u_i)\Big]
(1-v)+Mv,\cr
&(\theta_0,\ldots,\theta_{n+1})\in\Delta_{n+1},~~u_0,\ldots,u_{n+1}\in U,
{}~~v\in[0,1]\bigg\},\cr}\eqno(5.10)$$
where
$$\Delta_{n+1}=\Big\{(\theta_0,\ldots,\theta_{n+1});~~~\theta_i\in [0,1],
{}~~~\sum_{i=0}^{n+1}\theta_i=1\Big\}.$$

Since $F$ is continuous with compact convex values, the optimization
problem (5.4) subject to the boundary conditions (5.6) with dynamics
$$(\dot x,\dot x_0)\in F(t,x)\eqno(5.11)$$
admits an optimal solution.  The existence of a solution to the original
problem (5.1)-(5.3) will be proved by showing that the multifunction $F$
has the bang-bang property, for $(t,x,x_0)\in\Omega\times\R$.

To verify the concavity condition (C1), define the constant $\rho>0$ by
$${1\over\rho}=\sup\big\{ |\xi|;~~~\xi\in\partial_x\alpha(t,x),~~~(t,x)\in
\Omega\big\}.\eqno(5.12)$$
Given $(t,x)\in \Omega$, $(y,y_0)\in F(t,x)$, assume first $y_0<M-1$.
Choose any $\xi\in\partial_x\alpha(t,x)$ and consider the linear map
$$\Phi(z)\doteq\big(y+A(t)(z-x),~y_0+\xi\cdot(z-x)\big).
\eqno(5.13)$$
{}From (5.12) and the assumption on $y_0$ it follows that
$$y_0+\xi\cdot(z-x)\leq M\qquad~~\forall z\in\overline B(x,\rho).
\eqno(5.14)$$
Let $y,y_0$ be as in (5.10), for some $\theta_i, u_i,v$.  Using the
concavity of $\alpha$ we then obtain
$$y+A(t)(z-x)=A(t)z+\sum_{i=0}^{n+1}\theta_if(t,u_i),$$
$$\eqalign{y_0+\xi\cdot(z-x)&\geq\alpha(t,x)+\sum_{i=0}^{n+1}\theta_i
\beta(t,u_i)+\xi\cdot(z-x)\cr
&\geq\alpha(t,z)+\sum_{i=0}^{n+1}
\theta_i\beta(t,u_i).\cr}$$
This, together with (5.14), implies $\Phi(z)\in F(t,z)$.

Next, assume $y_0\in [M-1,~M]$.  Then the definition (5.8) implies
that the map
$$\Psi(z)\doteq\Big(y+A(t)(z-x),~y_0\Big)$$
is again a selection of $F(t,\cdot)$.
\v
We now turn to the condition (C2).  Fix $(t,x)\in \Omega$,~ $\tilde
y=(y,y_0)\in F(t,x)$,~ $\ve>0$, and choose $\xi\in\partial_x\alpha(t,x)$.
Write $\tilde y$ as a convex combination
$$\tilde y=\sum_{j=1}^\nu\vartheta_j\tilde y_j,\qquad~~~\vartheta_j\in(0,1],
{}~~~~\sum_{j=1}^\nu\vartheta_j=1,$$
with $\tilde y_j\in ext F(t,x)$, and define
$$\tilde y_j'=\tilde y_j+\ve '(\tilde y-\tilde y_j),$$
choosing $\ve '\in (0,1]$ such that
$$h\big(\tilde y_j',~F(t,x)\big)<\ve\qquad\quad\forall j.\eqno(5.15)$$
We now distinguish two cases.

If $y_0<M$, then $\tilde y_j'=(y_j',y_{0,j}')$ satisfies
$y_{0,j}'<M$ for all $j$.   Hence, for $|z-x|$ sufficiently small, the
maps
$$\Phi_j(z)\doteq\big(y_j'+A(t)(z-x),~y_{0,j}'+\xi\cdot(z-x)\big)$$
are affine selections of $F$.  Moreover, (5.15) and the upper semicontinuity
of $h$ imply
$$h\big(\Phi(z,z_0),~F(t,z)\big)<\ve$$
for all $z$ in a neighborhood of $x$.

On the other hand, if $y_0=M$, then $y_{0,j}=y_{0,j}'=M$ for all
$j$.  Hence the maps
$$\Phi_j(z)\doteq\big(y_j'+A(t)(z-x),~M\big)$$
are affine selections of $F$ and satisfy (5.16) for $|z-x|$ sufficiently
small.
\v
Applying Theorem 2, we now obtain the existence of an optimal solution $(x^
*,x_0^*):[0,T]\mapsto\R^n\times\R$ to (5.4), (5.6), (5.11), with the
additional property
$$(\dot x^*,~\dot x^*_0)(t)\in ext F(t,x^*(t))\qquad\hbox{for a.e.}
\quad t\in [0,T].$$
The representation (5.10) and the selection theorem [11] now
yield the existence of
some measurable $u^*:[0,T]\mapsto U$ such that
$$\left\{\eqalign{\dot x^*(t)&=A(t)x^*(t)+f(t,u^*(t)),\cr
x_0^*(t)&\in\big\{\alpha(t,x^*(t))+\beta(t,u^*(t)),~~M\big\}\cr}\right.$$
for almost every $t$.  Since the terminal value $x^*_0(T)$ is minimized,
by (5.8) we must have $\dot x^*_0(t)<M$ almost everywhere.
Therefore, $x^*$ is an optimal trajectory for the original system
(5.1), corresponding to the control $u^*$.
\vsk
\n {\medbf 6 - A counterexample.}
\v
The following example shows how the bang-bang property may fail, if some of
the assumptions in Theorem 2 or in Proposition 2 are not satisfied.
More general results concerning systems of this form can be found in [17, 18].
\v

On $\R^2$, consider the control system
$$\dot x(t)=f(x(t))+g(x(t))u(t), \qquad u(t)\in [-1,1],\eqno(6.1)$$
with
$$f(x_1,x_2)=(x_1,x_2),\qquad g(x_1,x_2)=(1,x_1).$$
For $t\in [0,1]$, the trajectory $t\mapsto (0,e^t)$, corresponding to the
control $u(t)\equiv 0$, steers the system from $p=(0,1)$ to $q=(0,e)$.

Defining the auxiliary function
$$V(x_1,x_2)=x_2-{x_1^2\over 2},$$
a straightforward computation yields
$${d\over dt}V(x(t))=V(x(t))-{x^2_1(t)\over 2}$$
for every solution of (6.1).  This implies
$$V(x(t))\leq e^tV(x(0)),$$
with equality holding if and only if $x_1(s)=0$ for all $s\in [0,t]$.
In particular, the control $u\equiv 0$ is the only one
which steers the system from $p$ to $q$ in minimum time.  Hence the
multifunction
$$F(x_1,x_2)=\big\{ (x_1+u,~x_2+x_1u);~~|u|\leq 1\big\}$$
does not have the bang-bang property.  Observe that
\v
\item{(i)} For each $y=f(x)+g(x)\omega\in F(x)$, defining
$$A=\left(\matrix{1 & 0 \cr \omega & 1\cr}\right)
\qquad c=\left(\matrix{\omega\cr 0\cr}\right),$$
one checks that the condition (C1) in Theorem 2 holds.
However, the condition (C2) here fails.
\v
\item{(ii)} The multifunction $extF(x)=\big\{f(x)+g(x),~f(x)-g(x)\big\}$
satisfies both (C1) and (C2) in Theorem 2, but its values are not convex.
\v
\item{(iii)} Each set $F(x)$ is a segment.  Moreover, $F$ admits the
representation
$$F(x)=\big\{y;~~w\cdot y\leq \psi_w(x)\doteq\max_{|u|\leq 1}~w\cdot\big(
f(x)+g(x)u
\big)~\big\}.$$
Since $f,g$ are linear, each $\psi_w$ is convex.  However, this
representation does not satisfy all assumptions in Proposition 2.
\vfill\eject
\centerline{\medbf References}
\vs
\item{[1]} T. S. Angell, Existence of optimal control without convexity and
a bang-bang theorem for linear Volterra equations, {\it J. Optim. Theory
Appl.}  {\bf 19} (1976), 63-79.
\v
\item{[2]} J.P. Aubin and A. Cellina, ``Differential Inclusions",
Springer-Verlag, Berlin/New York, 1984.
\v
\item{[3]} A. Bressan, The most likely path of a differential inclusion,
{\it J. Differential Equations}  {\bf 88} (1990), 155-174.
\v
\item{[4]} A. Bressan and G. Colombo, Generalized Baire category and
differential inclusions in Banach spaces,  {\it J. Differential
Equations}  {\bf 76} (1988), 135-158.
\v
\item{[5]} A. Cellina, On the differential inclusion $x'\in [-1,~+1]$,
{\it Atti Accad. Naz. Lincei Rend. Cl. Fis. Mat. Natur. Ser. VIII}
(1980), 1-6.
\v
\item{[6]} A. Cellina and G. Colombo, On a classical problem of the calculus
of variations without convexity assumptions, {\it Ann. Inst. Henri
Poincar\'e} {\bf 7} (1990), 97-106.
\v
\item{[7]} A. Cellina and F. Flores, Radially symmetric solutions of a class
of problems of the calculus of variations without convexity assumptions,
{\it Ann. Inst. Henri Poincer\'e}, to appear.
\v
\item{[8]} L. Cesari, ``Optimization - Theory and Applications",
Springer-Verlag, New York, 1983.
\v
\item{[9]} F. S. De Blasi and G. Pianigiani, Differential inclusions in
Banach spaces, {\it J. Differential Equations} {\bf 66} (1987),
208-229.
\v
\item{[10]} H. Hermes and J. P. LaSalle, ``Functional Analysis and Time
Optimal Control", Academic Press, 1969.
\v
\item{[11]} K. Kuratowski and C. Ryll-Nardzewski, A general theorem
on selectors, {\it Bull. Acad. Pol. Sc. Math. Astr. Phys.}
{\bf 13} (1965), 397-403.
\v
\item{[12]} L. W. Neustadt, The existence of optimal controls in the absence
of convexity conditions, {\it J. Math. Anal. Appl.}, {\bf 7} (1963),110-117.
\v
\item{[13]} J. P. Raymond, Th\`eor\'emes d'existence pour des probl\`emes
variationnels non convexes, {\it Proc. Royal Soc. Edimburgh} {\bf 107A}
(1987), 43-64.
\v
\item{[14]} J. P. Raymond, Existence theorems in optimal control problems
without convexity assumptions, {\it J. Optim. Theory Appl.} {\bf 67}
(1990), 109-132.
\v
\item{[15]} L. M. Sonneborn and F. S. Van Vleck, The bang-bang principle
 for linear control systems, SIAM {\it J. Control} {\bf 2} (1965)
, 151-159 .
\v
\item{[16]} M. B. Suryanarayana, Existence theorems for optimization problems
concerning linear hyperbolic partial differential equations without
convexity, {\it J. Optim. Theory Appl.} {\bf 19} (1976), 47-62.
\v
\item{[17]} H. J. Sussmann, The ``Bang Bang" problem for certain control
systems in $GL(n,\R)$, {\it SIAM J. Control}  {\bf 10} (1972), 470-476.
\v
\item{[18]} H. J. Sussmann, The structure of time-optimal trajectories for
single-input systems in the plane: the ${\cal C}^\infty$ nonsingular case,
{\it SIAM J. Control} {\bf 25} (1987), 433-465.
\v
\item{[19]} A. A. Tolstonogov, Extreme continuous selectors of multivalued
maps and the ``bang-bang" principle for evolution inclusions, {\it
Soviet Math. Dokl.}  {\bf 317} (1991).
\vs
\bye